\documentclass{elsarticle}
\usepackage{geometry,graphicx,amssymb,amsfonts,amsmath,bm,url}

\usepackage{graphicx}
\usepackage{algorithm}
\usepackage{algpseudocode}

\usepackage{latexsym,bm,amssymb,amsfonts,amsbsy,amsmath,graphicx,color,marvosym}
\usepackage{pifont}
\usepackage{url}
\usepackage{amsmath}
\usepackage{amsfonts}
\usepackage{amssymb}
\usepackage{fancyhdr}
\usepackage{pstricks,pst-node,pst-text,pst-3d}

\def\G1{\hbox{$\displaystyle{\mbox{\ding{172}}}$}}

\newtheorem{remark}{Remark}

\def \RR {{\mathbb{R}}}

\DeclareMathOperator{\MSE}{\overline{{E^2}}}
\DeclareMathOperator{\MSEj}{\overline{{E^2_j}}}
\DeclareMathOperator{\MSEzero}{\overline{{E^2_0}}}
\DeclareMathOperator{\MSEN}{\overline{{E^2_N}}}
\DeclareMathOperator{\diag}{diag}
\DeclareMathOperator*{\argmin}{arg\,min}

\newrgbcolor{marrone}{0.6156 0.3164 0.0313}
\newrgbcolor{ts}{0.50 0.598 0.781}
\newrgbcolor{verde}{0 0.4 0}
\newrgbcolor{LBlue}{0.58 0.65 0.9}
\newrgbcolor{Blue}{0.1 0.00 0.60}
\newrgbcolor{rosso}{0.9633    0.1641    0.0039}

\journal{Journal of Computational and Applied  Mathematics}

\begin{document}
\begin{frontmatter}
\title{An entropy-based approach for a robust least squares spline approximation}

\author[luigi]{Luigi Brugnano}
\ead{luigi.brugnano@unifi.it}
\author[mimmo]{Domenico Giordano}
\ead{dg.esa.retired@gmail.com}
\author[felixgiorgia]{Felice Iavernaro\corref{cor1}}
\ead{felice.iavernaro@uniba.it}
\cortext[cor1]{Corresponding author}
\author[felixgiorgia]{Giorgia Rubino}
\ead{g.rubino33@studenti.uniba.it}

\address[luigi]{Dipartimento di Matematica e Informatica ``U. Dini'', Universit\`a di Firenze, Italy}
\address[mimmo]{ESTEC (retired), European Space Agency, Noordwijk, The Netherlands}
\address[felixgiorgia]{Dipartimento di Matematica, Universit\`a degli Studi di Bari Aldo Moro, Italy}



\begin{abstract}
We consider the weighted least squares  spline approximation of a noisy dataset. By interpreting the weights as a probability distribution, we maximize the associated entropy subject to the constraint that  the  mean squared error is prescribed to a desired (small) value. Acting on this error yields a robust regression method  that automatically  detects and removes outliers from the data during the fitting procedure, by assigning them a very small weight. We discuss the use of both spline functions and spline curves. A number of numerical illustrations have been included to disclose the potentialities of the maximal-entropy approach in different application fields. 


\end{abstract}

\end{frontmatter}

\section{Introduction}
\label{intro}
With the advent of computer-aided modern technology, sheer volumes of data need to be preprocessed in order to make them suitable for the subsequent data-driven tasks they are intended for. Real data are often affected by various imperfections, including noise, poor sampling, missing values and outliers. The automatic  identification and removal of these inconsistencies has become of paramount importance during the preprocessing phase of data,  since they  may significantly affect the predictive accuracy and efficiency of models such as those based upon single and multivariate regression, as well as of pattern recognition procedures resulting from machine learning and deep learning processes \cite{ZhWuCh03, Te99, GaLaDz00, ZiFi18}.  

Identification of corrupted data also play a fundamental role in automatic anomaly  detection, meant as the appearance of events or observations which are inconsistent with the pattern underlying a given dataset. Anomaly detection has become increasingly important in many application areas ranging from statistics, cyber security,  medicine, event detection in sensor networks, financial fraud and machine learning \cite{ChBaKu09}.

Outliers may be thought of extreme values that deviate significantly from the trend defined by the majority of the data points, possibly due to errors or rare events, and that can consequently worsen the performance of many data analysis algorithms. Classical outlier detection methods often rely on specific assumptions about the data's distribution. However, in many real-world scenarios, estimating such a distribution beforehand can be challenging due to the data's dependence on various unknown or complex factors and the presence of highly noisy sources. This limitation becomes apparent in vast collections of time series data,  especially within environmental investigation. Such a topic has recently garnered extensive research attention, especially in understanding the correlation between climate changes and the increasing severity of natural disasters \cite{CaCaHo22, FaBiXiWaGu23}.

Extending the study addressed in \cite{GiIa21} for the polynomial case, the present paper introduces a robust regression technique for spline approximation of both univariate and multivariate time series, considering scenarios where observations exhibit varying degrees of reliability (see \cite{FaMaTa22} for a related study).
In statistics, robust regression tries to overcome the limitations of the ordinary least squares when its underlying assumptions are violated, for example, due to the presence of outliers \cite{An08,RoLe03,St16,YuYa17}. 

The proposed procedure tackles the challenges posed by outliers and noise by formulating a weighted least squares problem that leverages the statistical concept of entropy. To this end, we adopt the normalization condition that the weights sum to one, which allows us to interpret them as a probability distribution.

In more detail, to mitigate the negative influence of outliers and noise on the resulting approximating curve, the procedure maximizes the entropy $H$ associated with the weights distribution,  under the constraint that the resulting weighted mean squared error takes a prescribed value lower than the one corresponding to a uniform weights distribution. Such a value may be either provided by the user, on the basis of what he would expect in absence of corrupted data, or  automatically detected during the implementation of the procedure. 

To better elucidate the role played here by entropy, we quote Jaynes \cite[page 97]{Ja67}:
\begin{quote} \em
\dots the distribution that maximizes $H$, subject to constraints which represent whatever information we have, provides the most honest description of what we know. The pro\-ba\-bi\-lity is, by this process, spread out as widely as possible without contradicting the available information. 
\end{quote}

Translating Jaynes' words in our context,  we may stress that the proposed approach ensures that as many data points as possible carry non-negligible weights, which results in maximizing the inlier set while adhering to the mean squared error constraint. To achieve this, the strategy assigns smaller weights  to points that are more likely to be considered outliers, effectively minimizing their influence on defining the final shape of the approximating spline curve. It is important to note that this weighting task is seamlessly integrated into the fitting procedure, resulting in a unified methodology that eliminates the need for a preprocessing phase. Similarly to the RANSAC algorithm \cite{FiBo81}, the entropy-based approach proves particularly effective in handling situations where a substantial portion of the data is corrupted. However, unlike the RANSAC algorithm, it boasts the advantage of being deterministic in nature. Furthermore, by reinterpreting the weights as probabilities, we can readily justify the use of entropy as a mathematical tool for effectively handling corrupted data points.

The paper is structured as follows: In Section \ref{sec:2}, we review the fundamental concepts related to weighted least squares spline approximation and introduce the corresponding notations. Section \ref{sec:3} presents a formal definition of the approximation problem using the entropy-based tool and proposes a simple algorithm to obtain the optimal solution for the constrained optimization problem. To demonstrate the functionality of the entropy tool, a few numerical illustrations are provided in Section \ref{sec:4}. In Section \ref{sec:5}, three examples involving real-world data are considered. Finally, in Section \ref{sec:6}, we draw conclusions based on the findings.

\section{Background}
\label{sec:2}
Consider a parametrized sequence of points $\{(t_i,y_i)\}_{i=1}^{m}$, where $t=(t_1,\dots,t_m)^\top$ is a non-decreasing sequence of real pa\-ra\-me\-ters and $y_i\in\RR^s$ the corresponding data points. In the statistics parlance the  sequence $\{(t_i,y_i)\}$ is often referred to as a multivariate time series. As is usual in this context, we introduce a change of variable that normalizes the data in $[0,1] \times [0,1]^s$:\footnote{In the sequel, all the operations and functions evaluations involving vectors are meant componentwise. For example, for a given vector $z=(z_1,\dots,z_k)^\top$ and a function $g:\RR\rightarrow \RR$, we have  $g(z)=(g(z_1),\dots,g(z_k))^\top$.}
$$ 
t_i \rightarrow \frac{t_i-t_{\mathrm{min}}}{t_{\mathrm{max}} - t_{\mathrm{min}}},\qquad y_i \rightarrow \frac{y_i-y_{\mathrm{min}}}{y_{\mathrm{max}} - y_{\mathrm{min}}}.
$$
where 
$$
t_{\mathrm{min}} = \min_{1\le i \le m}{t_i}, \quad t_{\mathrm{max}} = \max_{1\le i \le m}{t_i}
$$
and, denoting by $y_i(j)$ the $j$th entry of the vector $y_i$,
$$
y_{\mathrm{min}}(j) = \min_{1\le i \le m}{y_i(j)}, \quad y_{\mathrm{max}}(j) = \max_{1\le i \le m}{y_i(j)}, \quad j=1,\dots,s.
$$
Of course, one can revert to the original coordinates by employing the inverse transformations. We wish to fit the given data set by means of a spline curve $f$ of degree $d$ expanded along a B-spline basis $\{B_{j,d}(x)\}_{j=1}^n$, namely
\begin{equation}
\label{spline}
f(x,c)=\sum_{j=1}^n c_j B_{j,d}(x).
\end{equation}
Here, $c=(c_1^\top,\dots,c_n^\top)^\top \in  \RR^{sn}$ is a set of $n$ control points, each of length $s$, and the B-splines $B_j(x)$ are defined on a non-decreasing sequence of  $(d+1)$-regular knots
\begin{equation}
\label{regular}
0=x_1=\dots=x_{d+1}<x_{d+2}\le \ldots  \le x_n<x_{n+1}=\dots=x_{n+d+1}=1,
\end{equation}
via the three-terms recursive relation\footnote{If a division by zero occurs, the related term is neglected.}
$$
B_{j,d}(x)=\frac{x-x_j}{x_{j+d}-x_j}B_{j,d-1}(x)+\frac{x_{j+d+1}-x}{x_{j+d+1}-x_{j+1}}B_{j+1,d-1}(x),
$$
with
$$
B_{j,0}(x)=\left\{\begin{array}{ll} 1, & \mbox{if } x_j\le x<x_{j+1},\\0, & \mbox{otherwise}. \end{array}\right. 
$$
Besides the conditions at the end points in (\ref{regular}),  the $(d+1)$-regularity of the knot vector  also imposes that $n\ge d+1$ and $x_j<x_{j+d+1}$, for $j=1,\dots,n$, which are relevant assumptions for the B-splines linear independence property \cite{LyMaSp18}. In the sequel, for sake of simplicity, we will omit the second subscript in $B_{j,d}(x)$.

Now, for a given vector $w=(w_1,\dots,w_m)^\top$ of (positive) weights satisfying the normalization condition
\begin{equation}
\label{normw}
\sum_{i=1}^m w_i=1,
\end{equation}
we consider the weighted mean squared error
\begin{equation}
\label{wmse}
\MSE = \sum_{i=1}^m w_i||f(t_i,c)-y_i||_2^2 
\end{equation}
as an estimate of the approximation accuracy of a  spline function $f(x,c)$ in the form (\ref{spline}) to the given data set. Denoting  by $I_s$ the identity matrix of dimension $s$ and introducing the generalized Vandermonde matrix
$$
A=\begin{pmatrix}
B_1(t_1) & \cdots & B_n(t_1) \\ 
\vdots   &        & \vdots \\
B_1(t_m) & \cdots & B_n(t_m) \\ 
\end{pmatrix} \in \RR^{m\times n},
$$
the vector $y=(y_1^\top,\dots,y_m^\top)^\top \in \RR^{sm}$ and the diagonal matrix $W=\diag(w_1,\dots,w_m)$,
(\ref{wmse}) may be cast in two equivalent forms that will be conveniently exploited for calculation and implementation purposes:
\begin{equation}
\label{wmse1}
\begin{array}{rcl}
\MSE &=& ( f(t,c)-y)^\top (W\otimes I_s) ( f(t,c)-y) \\[.15cm] 
         &=&\| (\sqrt{W}\otimes I_s)( f(t,c)-y)\|_2^2 \\[.15cm] 
         &=&\| (\sqrt{W}\otimes I_s)( (A\otimes I_s)c - y)\|_2^2\\[.15cm] 
\end{array} 
\end{equation}
and, denoting by $e_s=(1,\dots,1)^\top$  the unit vector of length $s$,
\begin{equation}
\label{wmse2}
\begin{array}{rcl}
\MSE &=& (w \otimes e_s)^\top (f(t,c)-y)^2 \\[.15cm] 
         &=& (w \otimes e_s)^\top ( (A\otimes I_s)c - y)^2.
\end{array} 
\end{equation}
For a prescribed choice of weights, the {\em least squares approximation problem} consists in finding the (vector) coefficients $c_j$ such that the corresponding weighted mean squared error (\ref{wmse1}) is minimized.  As is well known, differentiating (\ref{wmse1}) with respect to $c$, this requirement leads to the normal system  
\begin{equation}
\label{normal}
(A^\top W A \otimes I_s) c =(A^\top W \otimes I_s)y,
\end{equation}
which results from computing the stationary points of $\MSE$ regarded as a function of $c$.

Under the assumption that for any $j=1,\dots,n$ a $t_{i_j}$ exists such that $B(t_{i_j})\not = 0$,  matrix $A^\top W A$ is positive definite and the Cholesky factorization may be employed to transform (\ref{normal}) into a couple of triangular systems. More in general, also to prevent a worsening of the conditioning, one avoids the left multiplication by the matrix $A^\top$ and directly deals with the least squares solution of the overdetermined system 
\begin{equation}
\label{overdet}
(\sqrt{W}A \otimes I_s) c = \sqrt{W}y.
\end{equation} 
In such a case, application of the $QR$ factorization algorithm  with column pivoting,  or the SVD decomposition  to the rectangular matrix $\sqrt{W}A$ may be considered to solve the associated least squares problem.   
\begin{remark} \em
In the event that the components $y_i(j)$, $j=1,\dots,s$ are affected by sources of noise of different size depending on $j$, one could improve (\ref{wmse}) by allowing a different weight for each component of the error $f(t_i,c)-y_i$. This is tantamount to consider a vector of weights $w$ of length $ms$ and the related  mean squared error defined as 
\begin{equation}
\label{wmse_improved}
\MSE =  w^\top ( f(t,c)-y)^2 \equiv  \| \sqrt{W}( f(t,c)-y)\|_2^2,
\end{equation}
with $W=\diag(w)$. In the numerical tests discussed in Sections \ref{sec:4} and \ref{sec:5} both approaches showed pretty similar results, so we only included those relying on (\ref{wmse}).
\end{remark}

In the sequel, $\MSE_{\mathrm{uw}}$ will denote the  mean squared error  resulting from the ordinary least squared (OLS) approximation defined on the uniform  weights distribution  $w_i=1/m$, namely  
\begin{equation}
\label{Euw}
\MSE_{\mathrm{uw}} = \frac{1}{m}  \sum_{i=1}^m ( f(t_i,\bar c)-y_i)^2,
\end{equation}
where $\bar c$ satisfies the normal linear system (\ref{normal}) with $W=I_m/m$, $I_m$ being the identity matrix of dimension $m$.

\section{Maximum entropy weighted least squares spline approximation}
\label{sec:3}
The use of a weighted mean squared error is helpful when the data highlight different level of accuracies, due to the presence of noise and/or outliers. In such a case, it would be appropriate to attach large weights to very accurate data points  and small weights to data points which are most likely affected by a high level of inaccuracy. In fact, a weight $w_i$ approaching zero makes the corresponding data point $y_i$ irrelevant for the purpose of the fitting procedure. On the other hand, increasing the size of $w_i$ will make $f(t_i,c)$ closest to $y_i$. It turns out that, under the normalization condition (\ref{normw}), the WLS approximation will mimic the OLS one applied to the subset of data carrying relatively large weights.  

By exploiting an entropy-based argument, the  {\em maximum entropy weighted least squares} (MEWLS) approximation tries to devise an automatic, easy-to-understand and effective procedure for assigning the correct weight to each data point during the fitting procedure. The MEWLS approach based on splines approximating functions in the form (\ref{spline}) is defined by the following set of equations ($e_m$ stands for the unit vectors of length $m$):
\begin{equation}
\label{mels}
\begin{array}{rcl}
\mbox{maximize}     &~& -w^\top \log w, \\[.1cm]
\mbox{subject to:}  &~&  w^\top e_m=1, \\
                    &~&  (w \otimes e_s)^\top ( f(t,c)-y)^2 = \MSE.
\end{array}
\end{equation}
In other words, we wish to maximize the entropy function
\begin{equation}
\label{Hw}
H(w)=-w^\top \log w = -\sum_{i=1}^m w_i \log w_i
\end{equation}
associated with a weights distribution $w$ satisfying the normalization condition $\sum_i w_i=1$,  subject to the constraint that the corresponding  mean squared error attains a prescribed value $\MSE$. 

As is well known, problem (\ref{mels}),  deprived of the second constraint, admits the solution $w_i=1/m$, which leads us back to the ordinary least squares problem with uniform weights and associated means squared error  $\MSE_{\mathrm{uw}}$. Clearly, the very same solution is obtained when solving the complete set of equations in (\ref{mels}) under the choice  $\MSE=\MSE_{\mathrm{uw}}$, so (\ref{mels}) contains the ordinary least squares problem as a special instance. By setting  $\MSE$  to a  suitable value lower than the  mean squared error $\MSE_{\mathrm{uw}}$, the weights selection technique based upon the maximal-entropy argument epitomized by (\ref{mels})  is aimed at mitigating  the effect of outliers and noise in the data while solving the weighthed least squares problem. To highlight the relation between  $\MSE$ and $\MSE_{\mathrm{uw}}$, we assume in the sequel
\begin{equation}
\label{redfac}
\MSE = \frac{1}{r} \MSE_{\mathrm{uw}}
\end{equation}
where $r>1$ is a suitable reduction factor. 

According to the Lagrange multiplier theorem, we compute the stationary points of the Lagrangian function
\begin{equation}
\label{L}
{\mathcal {L}}(w,c,\lambda_1,\lambda_2) = w^\top \log w + \lambda_1 (w^\top e_m-1) 
                                         + \lambda_2\left((w \otimes e_s)^\top ( f(t,c)-y)^2 -\MSE\right). 
\end{equation}
Differentiating, we get:
\begin{eqnarray}
\displaystyle \frac{\partial {\mathcal {L}}}{\partial w} &=& e_m + \log w +\lambda_1 e_m \notag   \\
&& +\lambda_2 \left((I_m \otimes e_s)^\top(f(t,c)-y)^2\right), \label{dLw} \\[.3cm]
\displaystyle \frac{\partial {\mathcal {L}}}{\partial c} &=&  \displaystyle 2\lambda_2\left( (A^\top W A \otimes I_s) c - (A^\top W\otimes I_s)y\right), \label{dLc} \\[.3cm]
\displaystyle \frac{\partial {\mathcal {L}}}{\partial \lambda_1} &=&  w^\top e_m-1, \notag  \\[.3cm]
\displaystyle \frac{\partial {\mathcal {L}}}{\partial \lambda_2} &=&  (w \otimes e_s)^\top (f(t,c)-y)^2 -\MSE. \notag 
\end{eqnarray}
The last term in (\ref{dLw}) is the vector of length $m$
\begin{equation}
\label{last-term}
\lambda_2\left( ||f(t_1,c)-y_1||_2^2, \ldots, ||f(t_m,c)-y_m||_2^2\right)^\top,
\end{equation}
while (\ref{dLc}) comes from the equivalence of formulae (\ref{wmse1}) and (\ref{wmse2}), after observing that the first two terms in the Lagrangian (\ref{L}) do not depend on the spline coefficients $c_i$. The stationary points of $\cal L$ are the solutions of the following set of $n+m+2$ equations in as many unknowns $c\in\RR^n$, $w\in \RR^m$, $\lambda_1$ and $\lambda_2$:

\begin{eqnarray}
\displaystyle (A^\top W A \otimes I_s) c - (A^\top W\otimes I_s)y &=& 0, \label{normal-eq} \\[.3cm]
\displaystyle   (w \otimes e_s)^\top (f(t,c)-y)^2 -\MSE &=& 0,  \label{mse-eq} \\[.3cm]  
\displaystyle e_m + \log w +\lambda_1 e_m + \lambda_2\left((I_m \otimes e_s)^\top(f(t,c)-y)^2\right)&=&0, \label{weights-eq} \\[.3cm]
\displaystyle   w^\top e_m-1 &=&0. \label{normalization-eq} 
\end{eqnarray}
By exploiting the weights normalization condition (\ref{normalization-eq}), we can easily remove the unknown $\lambda_1$. To this end, we first recast equation (\ref{weights-eq}) as 
$$
w=\exp( -(1+\lambda_1)) \cdot \exp\left(-\lambda_2\left((I_m \otimes e_s)^\top(f(t,c)-y)^2\right)\right).
$$
Multiplying both sides by $e_m^\top$ and taking into account (\ref{normalization-eq}) and (\ref{last-term}) yields
$$
1=\exp( -(1+\lambda_1)) \cdot Q(c,\lambda_2), \qquad \mbox{with } Q(c,\lambda_2)=\sum_{i=1}^m \exp\left( -\lambda_2||f(t_i,c)-y_i||_2^2\right)
$$
and hence
\begin{equation}
\label{weight1}
w=\frac{1}{Q(c,\lambda_2)}\cdot \exp\left(-\lambda_2\left((I_m \otimes e_s)^\top(f(t,c)-y)^2\right)\right)
\end{equation}
that will replace (\ref{weights-eq}) and (\ref{normalization-eq}). Plugging (\ref{weight1}) into (\ref{mse-eq}) we  arrive at the final shape of the system to be solved:
\begin{eqnarray}
\displaystyle (A^\top W A \otimes I_s) c - (A^\top W\otimes I_s)y &=& 0, \label{normal_eq} \\[.3cm]
\displaystyle   \sum_{i=1}^m ||f(t_i,c)-y_i||_2^2 \cdot \exp\left( -\lambda_2||f(t_i,c)-y_i||_2^2\right)  - \sum_{i=1}^m \exp\left(-\lambda_2||f(t_i,c)-y_i||_2^2\right)\MSE &=& 0,  \label{mse_eq} \\[.3cm]  
\displaystyle w-\frac{1}{Q(c,\lambda_2)}\cdot \exp\left(-\lambda_2\left((I_m \otimes e_s)^\top(f(t,c)-y)^2\right)\right)&=&0. \label{weights_eq} 
\end{eqnarray}
Before facing the question of how to solve the system numerically,  a few remarks are in order: 
\begin{itemize}
\item (\ref{normal_eq}) is nothing but the normal linear system one would get when handling the least squares problem with constant weights (see (\ref{normal})). It can be therefore expressed as the overdetermined system (\ref{overdet}) which has to be solved in the least squares sense;
\item (\ref{mse_eq}) is a scalar equation that, for a given vector $c$, may be easily solved with respect to the Lagrange multiplier $\lambda_2$ via a Newton or Newton-like iteration; 
\item equation (\ref{weights_eq}) is explicit with respect to the unknown $w$, for given $\lambda_2$ and $c$.  
\end{itemize}
Therefore, a quite natural technique to solve the nonlinear system (\ref{normal_eq})--(\ref{weights_eq}) is yielded  by the  hybrid iteration summarized in Algorithm \ref{Alg1} ($tol$ is an input tolerance for the stopping criterion).

\begin{algorithm}
\begin{algorithmic}[1]
\State initially, set  $W^{(0)} \gets\frac{1}{m}I_m$, $\lambda_2^{(0)}\gets 0$; $k\gets0$; 
\Repeat:
\State {$k \gets k+1$;}
\State $\displaystyle  c^{(k)} \gets \argmin_c{||\sqrt{W^{(k-1)}}Ac-\sqrt{W^{(k-1)}}y||_2}$;
\State employ a Newton iteration scheme with initial guess $c^{(k)},\lambda_2^{(k-1)}$, to solve  (\ref{mse_eq}) and get $\lambda_2^{(k)}$;  \vspace*{.2cm}
\State $\displaystyle w^{(k)} \gets \frac{1}{Q(c^{(k)},\lambda_2^{(k)})}\cdot \exp\left(-\lambda_2^{(k)}\left((I_m \otimes e_s)^\top(f(t,c^{(k)})-y)^2\right)\right)$; \vspace*{.2cm}
\State $W^{(k)}=\diag(w^{(k)})$;
\Until $(|c^{(k)}-c^{(k-1)}|<tol) \quad \& \quad  (|\lambda_2^{(k)} - \lambda_2^{(k-1)}|<tol) \quad \& \quad  (|w^{(k)} -w^{(k-1)}|<tol) $;
\end{algorithmic}
\caption{Numerical procedure for  solving system (\ref{normal_eq})--(\ref{mse_eq})}
\label{Alg1}
\end{algorithm}
In order to improve the convergence properties of the nonlinear scheme,  we employ a continuation technique on $ \MSE$. In more detail, we define a sequence of increasing reduction factors 
$$
1=r_0<r_1<r_2<\dots<r_N =\frac{\MSE_{\mathrm{uw}}}{\MSE}
$$ 
and the corresponding sequence of mean squared errors   
\begin{equation}
\label{continuation}
\MSEj= \frac{1}{r_j}\MSE_{\mathrm{uw}}, \quad j=0,\dots,N,
\end{equation}
so that $\MSEzero=\MSE_{\mathrm{uw}}$ and $\MSEN=\MSE$. Then, for $j=0,\dots,N$, we perform lines 2--5 of Algorithm \ref{Alg1} taking care that the output quantities $c^{(k)}$, $\lambda_2^{(k)}$, $W^{(k)}$ obtained at step $j$ are used as input parameters for the subsequent step $j+1$. 

A further relevant motivation for employing such a continuation technique is that it generates a discrete family of homotopic curves, parametrized by $\MSEj$,  admitting the OLS and the MEWLS solutions as initial and final configurations respectively. Each element in this family brings a specific weights distribution (and entropy value) and acts as a starting guess for the subsequent approximation curve. Therefore, the overall procedure can be interpreted as an improvement on the OLS approximation in that,  by reducing the mean squared error progressively, it smoothly deforms the initial shape of the spline curve  to get rid of outliers.  An illustration is provided in the first example of the next section.    
   
Finally, it is worth noticing that the resulting weights may be exploited for classification purposes. Indeed, the original  data set $D$ may be split in two disjoint subsets: $D= D_1 \cup D_2$, where $D_1$ contains the inliers while $D_2$ identifies the outliers. To this end, given a small enough tolerance $tol$, one can set, for example, 
\begin{equation}
\label{D1D2}
D_2= \{(x_i,y_i)\in D ~|~ w_i < tol \cdot \max_j{w_j}\}, \qquad D_1=D - D_2.
\end{equation}

\section{Numerical illustrations}
\label{sec:4}
To showcase the potential of the MEWLS spline approximation, we present three numerical experiments using synthetic data points. The first experiment focuses on a spline function fitting problem, aiming at elucidating the continuation technique (\ref{continuation}) and the use of (\ref{D1D2}) for the automatic detection of outliers. The second and third examples involve approximating a set of data points with a spline curve in the plane and in 3D space, respectively. All the numerical tests have been implemented in Matlab (R2023a) on a 3.6 GHz Intel I9 core computer with 32 GB of memory. References to colors have been included for the online version of the manuscript.

\subsection{Example 1}
\label{subsec:4.1}
\begin{figure}[htb]
\begin{center}
\includegraphics[width=0.5\textwidth]{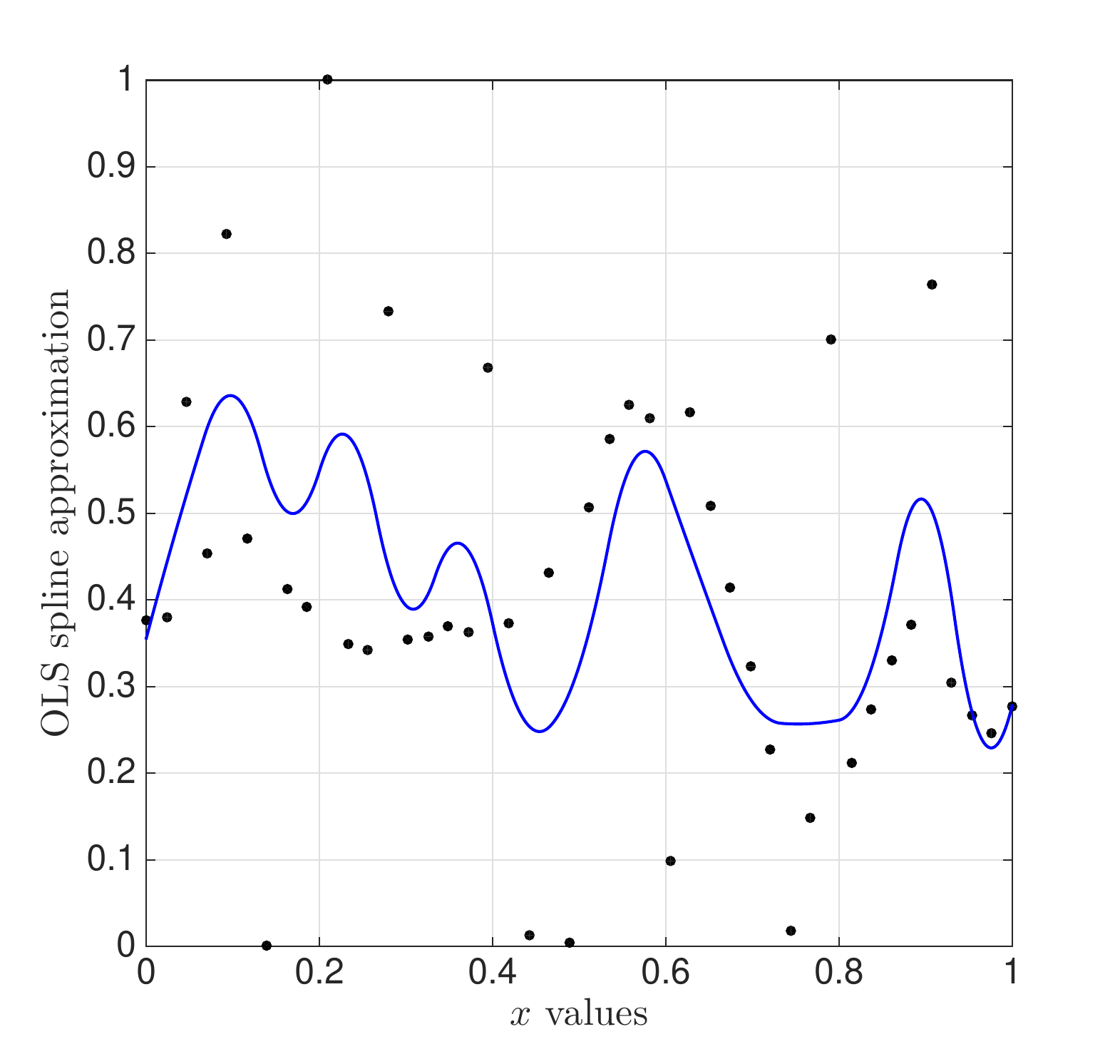}\includegraphics[width=0.5\textwidth]{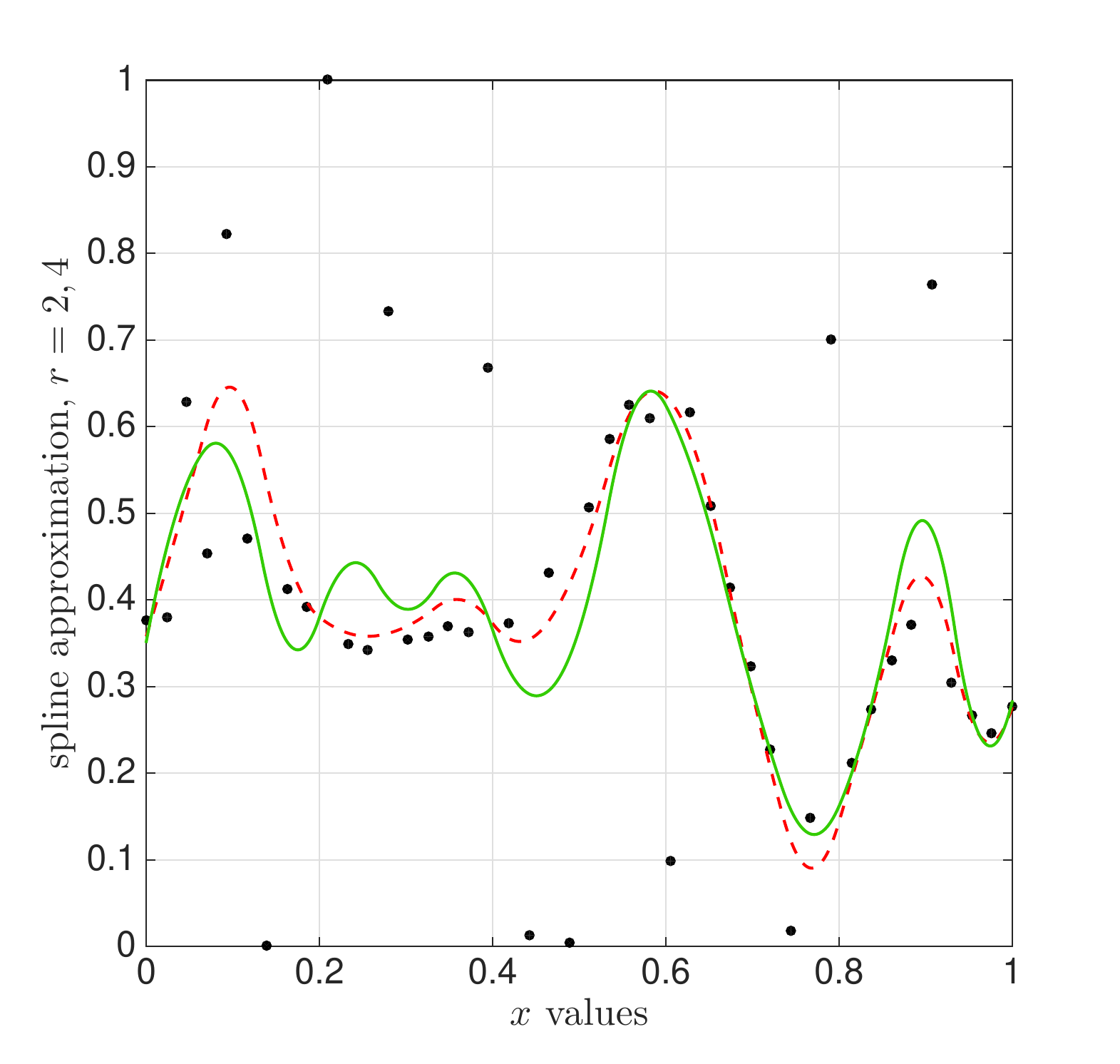} \\
\includegraphics[width=0.5\textwidth]{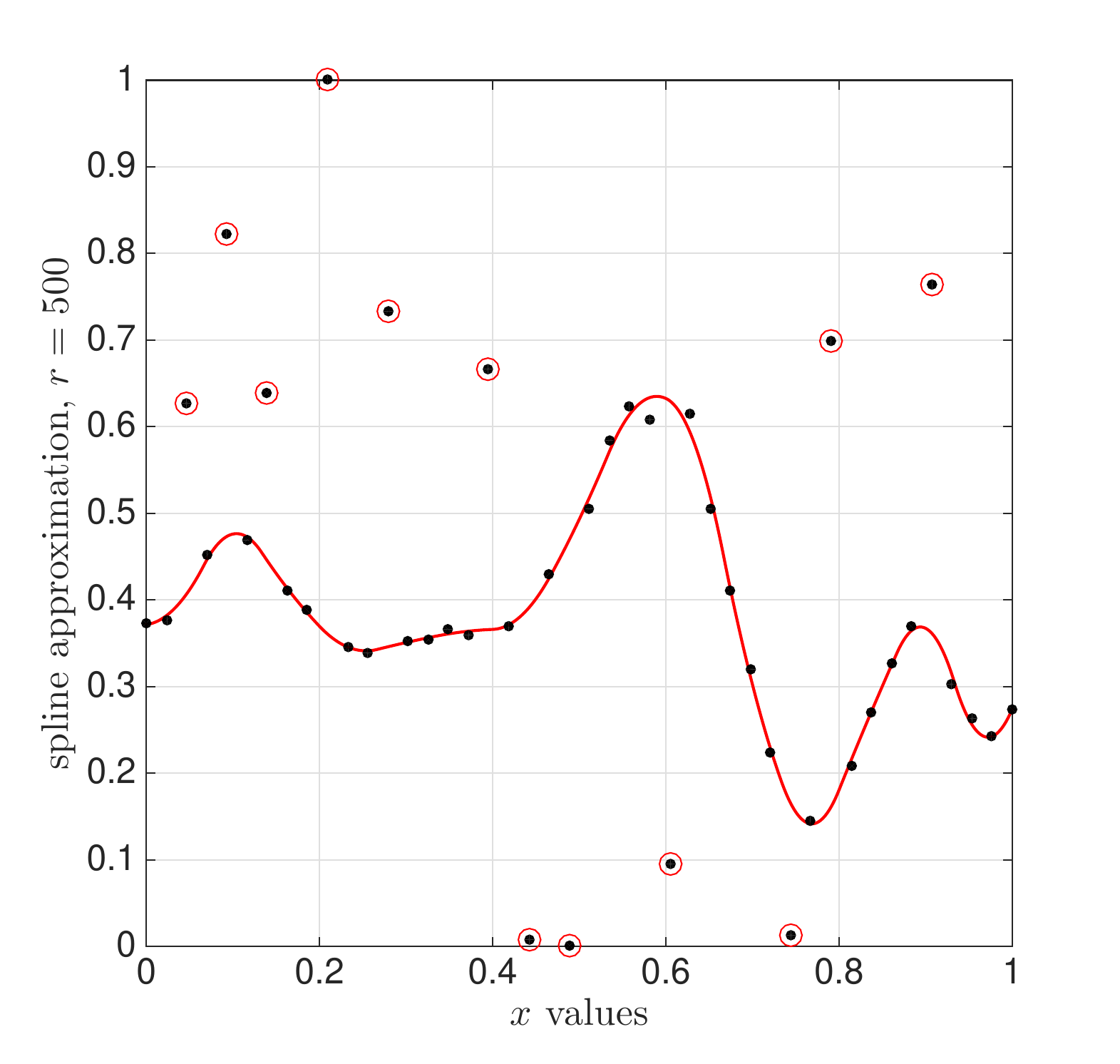}\includegraphics[width=0.5\textwidth]{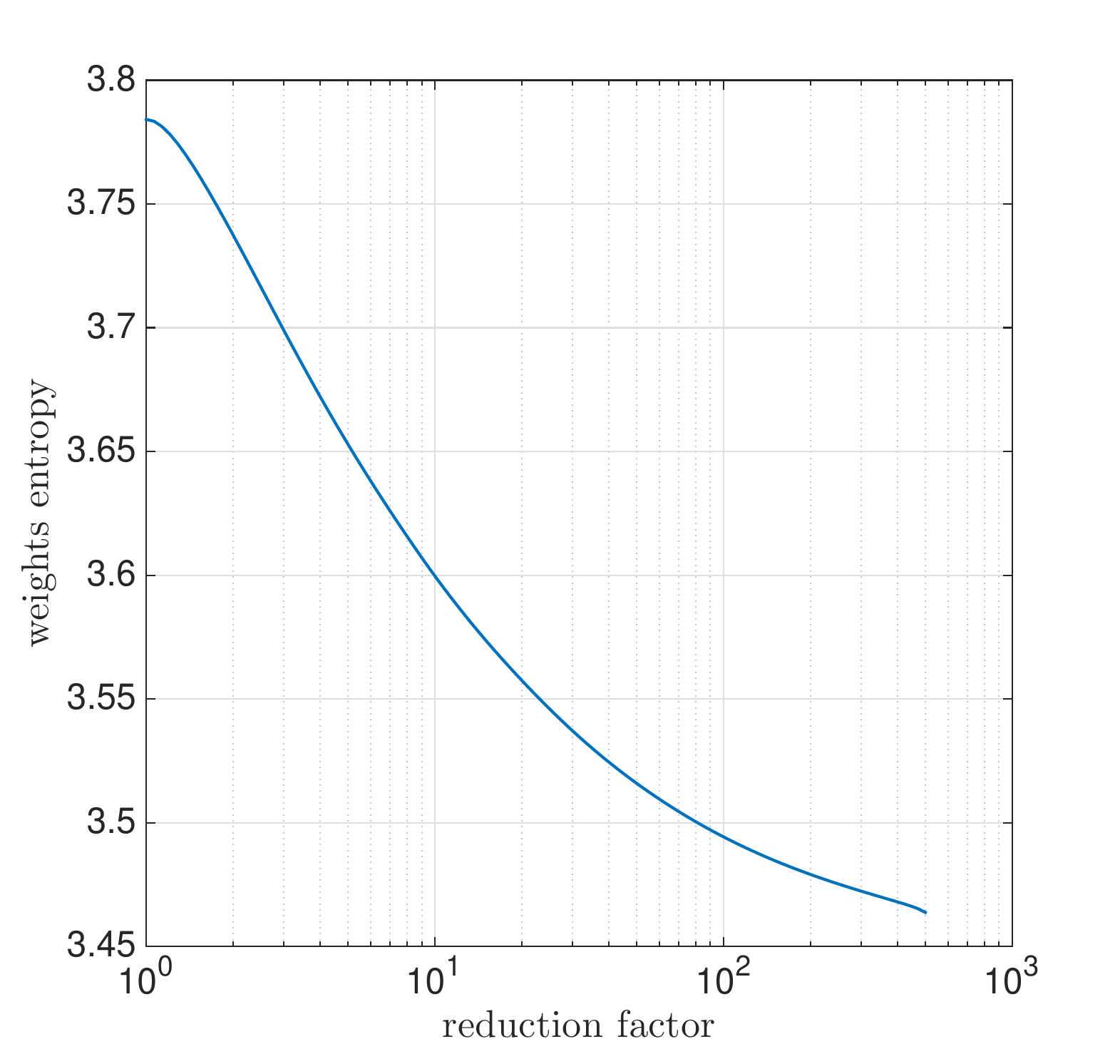}
\end{center}
\caption{Results obtained for Example 1. Top-left picture: a noisy data set revealing a pattern (dots) and its OLS spline approximation (blue line). Top-right picture: two homotopic spline functions corresponding to a reduction factor $r=2$ (solid green line) and $r=4$ (dashed red line). Bottom-left picture: final MEWLS spline approximation (red line) obtained by reducing the mean squared error by a factor $r=500$. Detected outliers, identified automatically through the use of formula (\ref{D1D2}), are indicated by dots surrounded with small circles. Bottom-right picture: the entropy associated with the distribution of weights, showcased as a function of the scaling factor $r$.}
\label{Fig_spline}
\end{figure}
We consider a dataset comprising $44$ points in the square $[0,1]\times [0,1]$, out of which 32 closely follow a given profile, while the remaining 12 consistently deviate from it. To fit the data, we employ a spline of degree $d=2$, defined on a regular and uniform knot sequence consisting of $20$ nodes, covering the interval $[0,1]$.

In the top-left picture of Figure \ref{Fig_spline}, we observe the data set along with the ordinary least squares approximation. We see that the OLS spline approximation fails to accurately reproduce the correct profile due to the strong influence of the 12 anomalous points. Therefore, we aim to improve the approximation by decreasing the weighted mean squared error while utilizing the maximal-entropy argument to make an optimal weights selection. To this end, we consider  a sequence of reduction factors distributed over the interval $[1,500]$. For graphical clarity, we set $N=50$ in (\ref{continuation}) to mimic the behavior of formula (\ref{redfac}), where the variable $r$ continuously varies within the specified interval. Algorithm \ref{Alg1} generates a sequence of $50$ homotopic functions with parameter $r\in[1,500]$.

The top-right picture of Figure \ref{Fig_spline} displays two such functions, one corresponding to $r=2$ ($\MSE= \MSE_{\mathrm{uw}}/2$, solid line), and the other to $r=4$ ($\MSE= \MSE_{\mathrm{uw}}/4$, dashed-line). As the reduction factor $r$ increases, the maximum entropy principle deforms the shape of the original OLS solution by adjusting the weights to ensure that the maximum number of points contribute while still adhering to the mean squared error constraint. 

In the bottom-left picture of Figure \ref{Fig_spline} we can see the final shape of the approximating spline, corresponding to $r=500$ ($\MSE= \MSE_{\mathrm{uw}}/500$). We can see that it nicely conforms to the profile underlying the given data set. The use of formula (\ref{D1D2}) with $tol=10^{-4}$ correctly detects $12$ outliers which are surrounded by small circles in the picture. 

Finally, the bottom-right picture of Figure \ref{Fig_spline} illustrates the behavior of the entropy (\ref{Hw}) as a function of the scaling factor $r$. As expected, reducing $\MSE$ results in a decrease of the entropy associated with the weights distribution. The appropriate choice of $\MSE$ depends on the context and, in particular, on the expected accuracy of the model in the absence of outliers. An automatic identification of a suitable value for $\MSE$ may be inferred by examining the rate of change in the spline approximations as the scaling factor $r$ increases, which is closely related to the behavior of the entropy function $H$ as a function of $r$. This aspect will be the subject of future research.

\subsection{Example 2}
\label{subsec:4.2}
We address the problem of approximating the arithmetic spiral defined by the equations 
$$
\left\{
\begin{array}{rcl}
x(t) &=& (a+bt)\cos(t),\\
y(t) &=& (a+bt)\sin(t), 
\end{array}
\right.
$$
with $a=1$, $b=4$, $t\in[-a/b,4\pi]$, which ensures that the spiral originates at the origin.  To this end, we create a data set consisting of $N=200$ points sampled along the spiral and then introduce random noise to  $100$ of them, specifically targeting the odd-numbered ones. In more detail, after setting $h=4\pi/(N-1)$, our data set is defined as follows:
$$
\left\{
\begin{array}{rcll}
t_i  &=& (i-1)h, \qquad &i=1,\dots,N, \\
(x_i,y_i)  &=& (x(t_i),y(t_i)), \qquad & \mbox{if $i$ is even}, \\
(x_i,y_i)  &=& (x(t_i)+\delta_x^{(i)},y(t_i)+\delta_y^{(i)}), \qquad & \mbox{if $i$ is odd},
\end{array}
\right.
$$
where $\delta_x^{(i)},\delta_y^{(i)} \in {\cal N} (0, \sigma^2)$ are random variables distributed normally with mean 0 and variance $\sigma^2=30$. Since, for the specified range of $t$, the spiral is entirely enclosed in the square $S=[-60,60]^2$, for  visualization clarity,  we iterate the generation of values  $\delta_x^{(i)},\delta_y^{(i)}$ until $(x_i,y_i)$ falls within $S$, for each odd index $i$. 

\begin{figure}[htb]
\begin{center}
\includegraphics[width=0.5\textwidth]{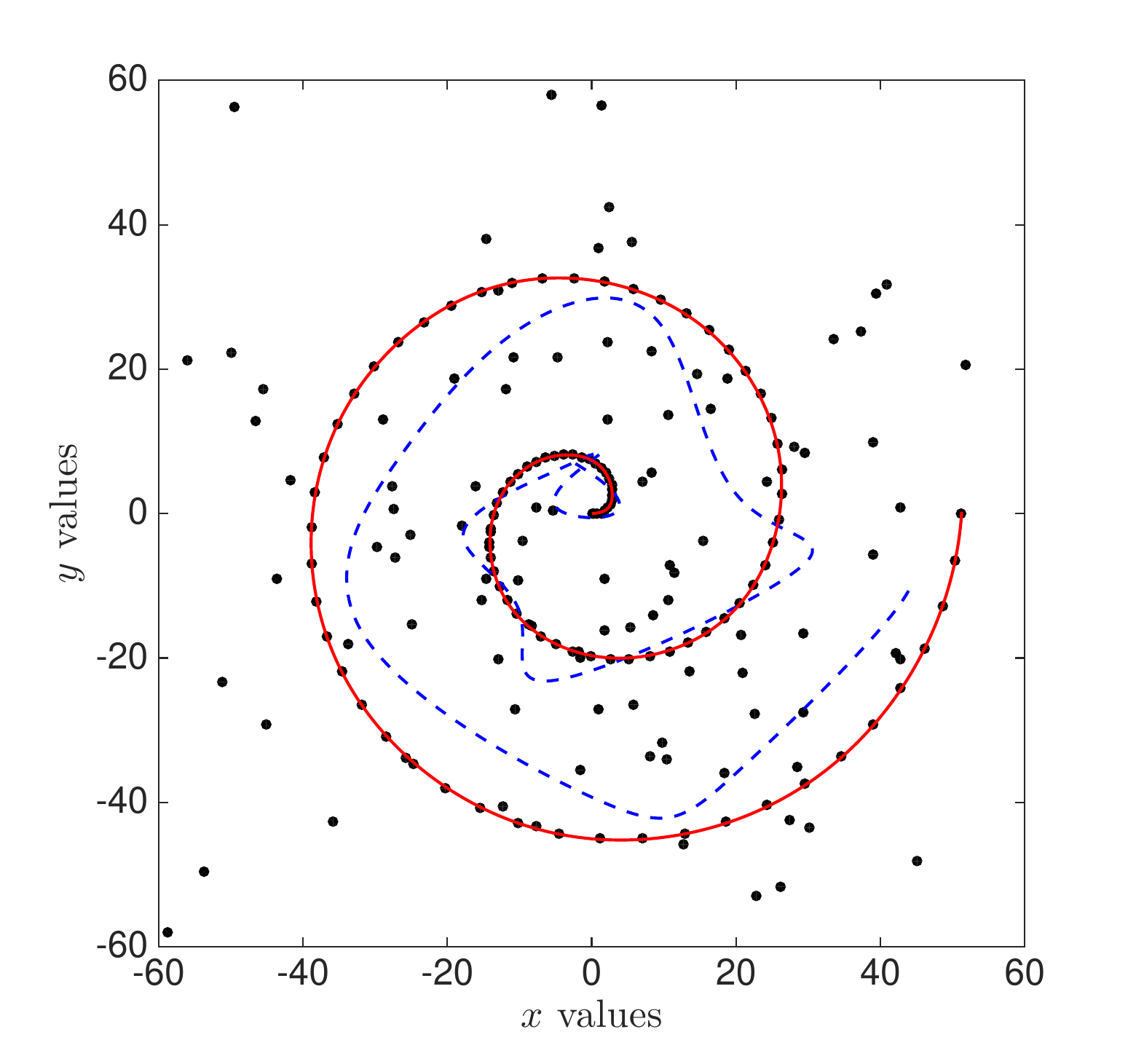}\includegraphics[width=0.5\textwidth]{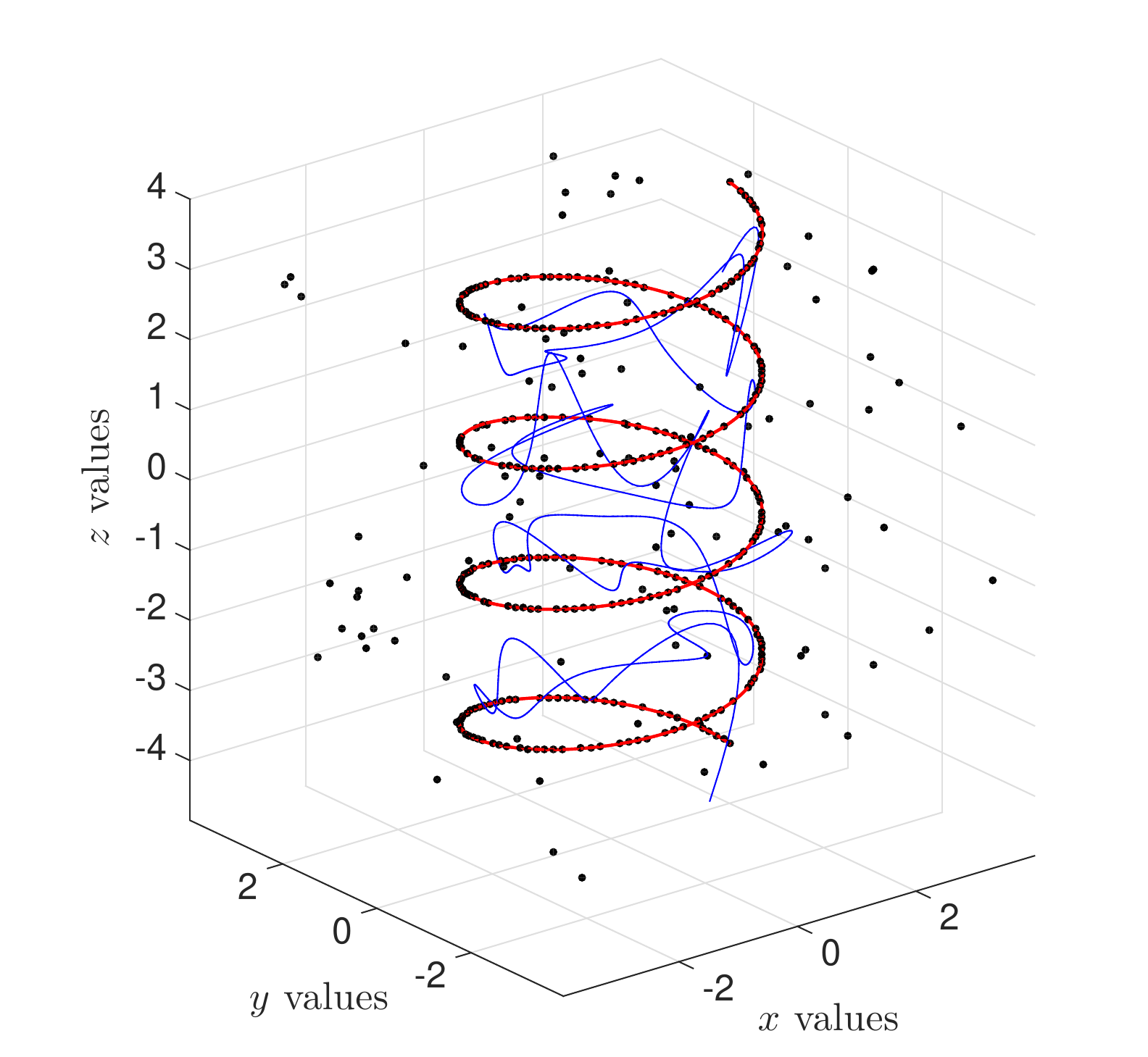} 
\end{center}
\caption{Results obtained for Examples 2 and 3. Left picture: a dataset comprising $200$ points, with half of them precisely aligned on an Archimedean spiral and the remainder introducing noise. Both OLS (dashed blue line) and MEWLS (solid red line) spline approximations are illustrated. Right picture: the data set consists of $400$ points, with $300$ of them following a circular helix pattern, while the remaining $100$ contribute as noise.  Both  OLS (irregular blue line) and MEWLS (red line) spline approximations are displayed.}
\label{Fig_spiral-helix}
\end{figure}

The left picture of Figure \ref{Fig_spiral-helix} portrays the dataset ${(x_i, y_i)}_{i=1}^N$ along with the spline approximations using ordinary least squares (dashed line) and maximum entropy weighted least squares  (solid line). Notably, while the OLS approximation struggles to capture the true spiral due to the presence of outliers, the MEWLS spline curve faithfully reproduces the unperturbed spiral $(x(t), y(t))$. 

\subsection{Example 3}
\label{subsec:4.3}
We replicate a procedure akin to the one executed in the prior spiral example but address our attention to a circular helix defined by the equations
$$
\left\{
\begin{array}{rcl}
x(t) &=& r\cos(2 \pi t),\\
y(t) &=& r\sin(2 \pi t), \\
z(t) &=& ct, 
\end{array}
\right.
$$
with $r=2, c=1$ and $t\in[-4,4]$, so the helix is enclosed in the cube $[-4,4]^3$. We begin with a data set ${(x_i, y_i, z_i)}_{i=1}^N$ consisting of $N=400$ points sampled along the helix but, differently to what was done in Example \ref{subsec:4.2}, we now introduce a random noise to a randomly chosen subset of these points. More precisely, we first compute a subset $\Omega$ obtained by randomly extracting $M$ points from the set of indices $\{1,2,\dots,N\}$. Then we define    
$$
\left\{
\begin{array}{rcll}
t_i  &=& (i-1)h, \qquad &i=1,\dots,N, \\
(x_i,y_i,z_i)  &=& (x(t_i),y(t_i),z(t_i)), \qquad & \mbox{if $i \not \in \Omega$}, \\
(x_i,y_i,z_i)  &=& (x(t_i)+\delta_x^{(i)},y(t_i)+\delta_y^{(i)},z(t_i)+\delta_z^{(i)}), \qquad & \mbox{if $i \in \Omega$}.
\end{array}
\right.
$$
Here, $\delta_x^{(i)}, \delta_y^{(i)}, \delta_z^{(i)} \in \mathcal{N} (0, 20)$ represent random variables drawn from a normal distribution with mean $0$ and variance $\sigma^2 = 20$. Again, for  visualization clarity,   for each odd index $i$ we iterate the generation of the perturbation values  $\delta_x^{(i)},\delta_y^{(i)}, \delta_z^{(i)}$ until $(x_i,y_i,z_i)$ falls within the cube $S=[-4,4]^3$. The right picture of Figure \ref{Fig_spiral-helix} displays the dataset ${(x_i, y_i, z_i)}_{i=1}^N$ along with the spline approximations using ordinary least squares (irregular solid line) and maximum entropy weighted least squares  (helix-shaped solid line). Again the MEWLS spline curve faithfully reproduces the shape of the original helix.  The results obtained in both this example and the previous one  underscore the effectiveness of MEWLS in successfully detecting and eliminating outliers from highly noisy datasets. Further instances based on real data are illustrated in the next section.

\section{A few applications to real data}
\label{sec:5}
\subsection{Approximating the main sequence in a Hertzsprung--Russell diagram}
\label{sec:5.1}
The Hertzsprung-Russell (HR) diagram is a graphical representation of stars, mapping the correlation between their absolute magnitudes or luminosities versus their color indices or temperatures, allowing astronomers to discern distinct patterns in stellar evolution \cite{Maoz16,CaOs17}. 

The absolute magnitude of a star is a measure of its intrinsic brightness or luminosity, unaltered by its distance from Earth. It is the apparent magnitude (brightness as seen from Earth) that a star would have if it were located at a standard distance of $10$ parsecs (about $32.6$ light-years) away. Essentially, the absolute magnitude allows astronomers to compare the luminosities of stars irrespective of their varying distances from us.

The B-V color index is a parameter that characterizes a star's color and temperature. It is the difference between the star's apparent magnitudes in the blue (B) and visual (V) parts of the electromagnetic spectrum. Blue stars have negative B-V values, while redder stars have positive values. This index is crucial in categorizing stars by their spectral types, indicating whether a star is hotter (blue) or cooler (red).

\begin{figure}[htb]
\begin{center}
\includegraphics[width=0.5\textwidth]{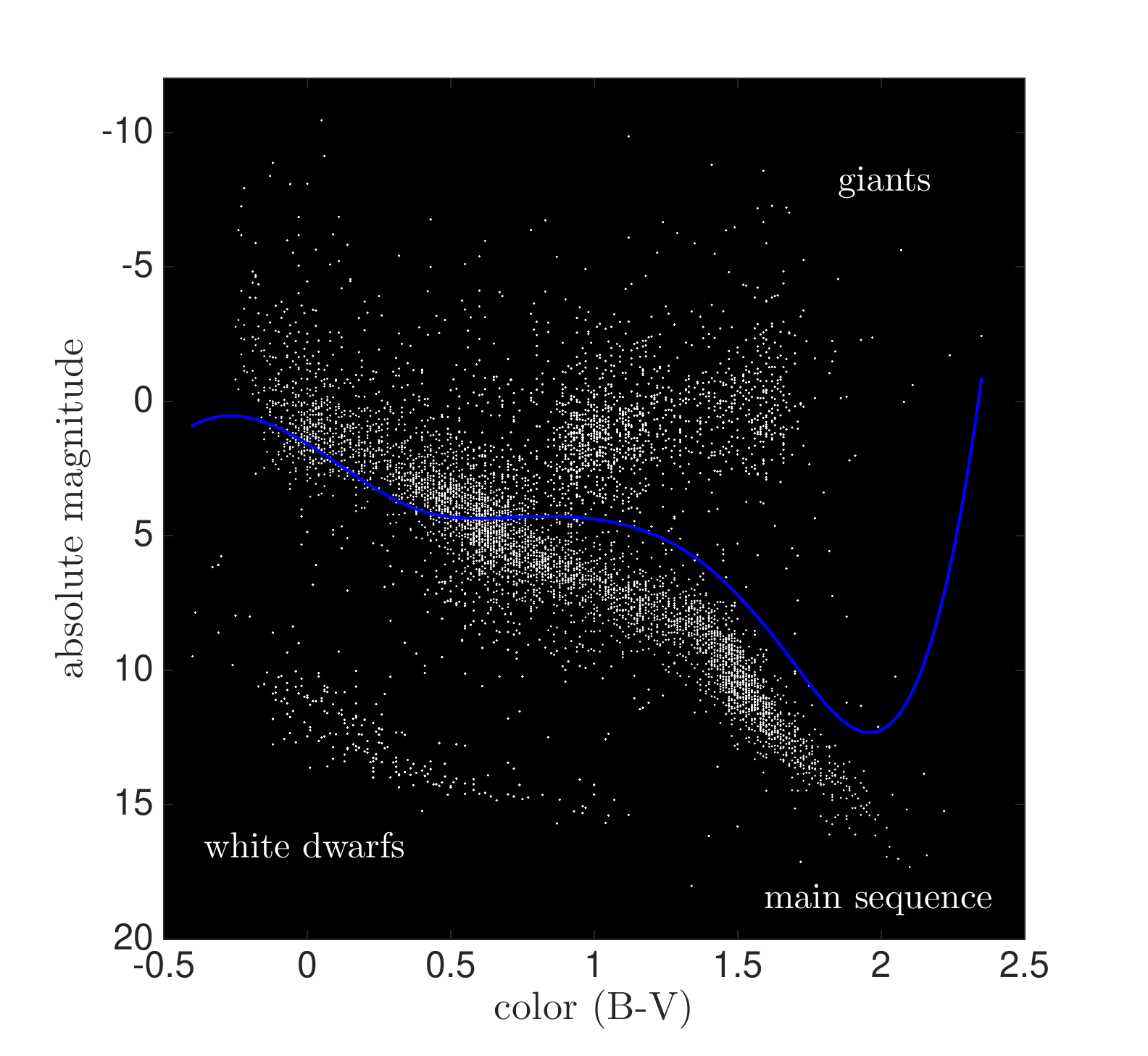}\includegraphics[width=0.5\textwidth]{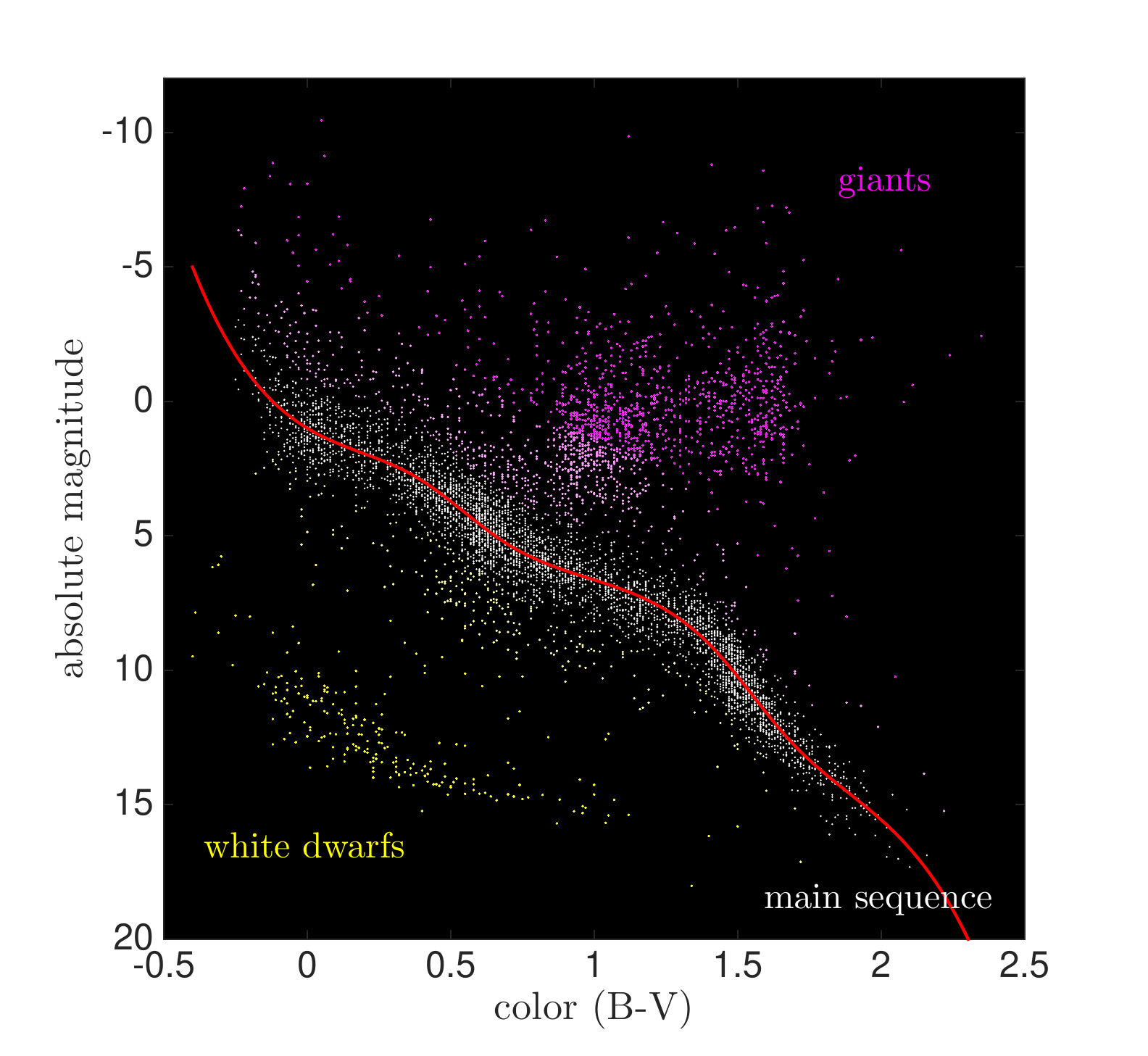} 
\end{center}
\caption{Hertsprung-Russel diagrams of the Yale dataset.  Left picture: ordinary least squares spline approximation (blue line).  Right picture: maximal-entropy least squares spline approximation (red line). The intensity of magenta and yellow colors is inversely proportional to the weight associated with each data point.}
\label{Fig_stars}
\end{figure}

Together, the absolute magnitude and B-V color index are vital tools in understanding stars' properties, evolutionary stages, and positions within the Hertzsprung-Russell diagram. As an example, the left picture of Figure \ref{Fig_stars} shows the HR diagram for the Yale Trigonometric Parallax Dataset \cite{AlLeHo95} comprising more than 6000 catalogued stars.
This astronomical resource  provides  measurements of stellar distances using the trigonometric parallax method, a technique employed to determine the distance to a star by measuring its apparent shift in position against more distant background stars as the Earth orbits the Sun. Besides observed parallaxes (in arcsec), the Yale catalogue also includes the B-V color index and the apparent V magnitude. The absolute magnitude is then obtained by means of the formula
$$
\mbox{absolute magnitude} = \mbox{apparent V magnitude} +5(\log_{10}(\mbox{observed parallax})+1).
$$

At its core, the diagram features a  continuous and well-defined band known as the main sequence. This band comprises the vast majority of genuine stars in the cosmos, including our own Sun with an absolute magnitude of 4.8 and a B-V color index of 0.66.

Located in the lower-left portion of the diagram are the white dwarfs, while the upper part accommodates the subgiants, giants, and supergiants. This layout visually captures the diverse stages of stellar evolution with the white dwarfs representing stars in their final stages of evolution.

One of the diagram's remarkable applications is in determining the distance between Earth and distant celestial objects like star clusters or galaxies. 

In this example, our aim is to accurately approximate the main sequence's shape using an appropriate spline curve and further categorize stars through color assignments. To achieve this, we employed a spline of degree $d=3$ along with a regular knot sequence 
$$
t=[0, 0, 0, 0, 0.286, 0.397, 0.658, 0.757, 1, 1, 1, 1].
$$
The left picture in Figure \ref{Fig_stars} displays the outcome of the ordinary least squares approximation (indicated by the blue line in the color image). This method evidently fails to accurately replicate the main sequence's distinctive form due to the presence of giants and white dwarfs. Conversely, the maximal-entropy least squares approximation successfully captures the main sequence's true shape. By determining the distribution of weights based on the entropy-driven procedure, we assigned distinct color gradients to each star. This color differentiation effectively highlights the discrepancies between these stars and those belonging to the main sequence. As the corresponding weights decrease, the intensity of magenta and yellow pixels progressively intensifies. This approach not only improves the accuracy of the main sequence representation but also facilitates the identification of stars that deviate from its expected characteristics.

\subsection{Detecting train rails in a railway infrastructure and surrounding environment}
\label{sec:5.2}

\begin{figure}[htb]
\begin{center}
\includegraphics[width=0.93\textwidth]{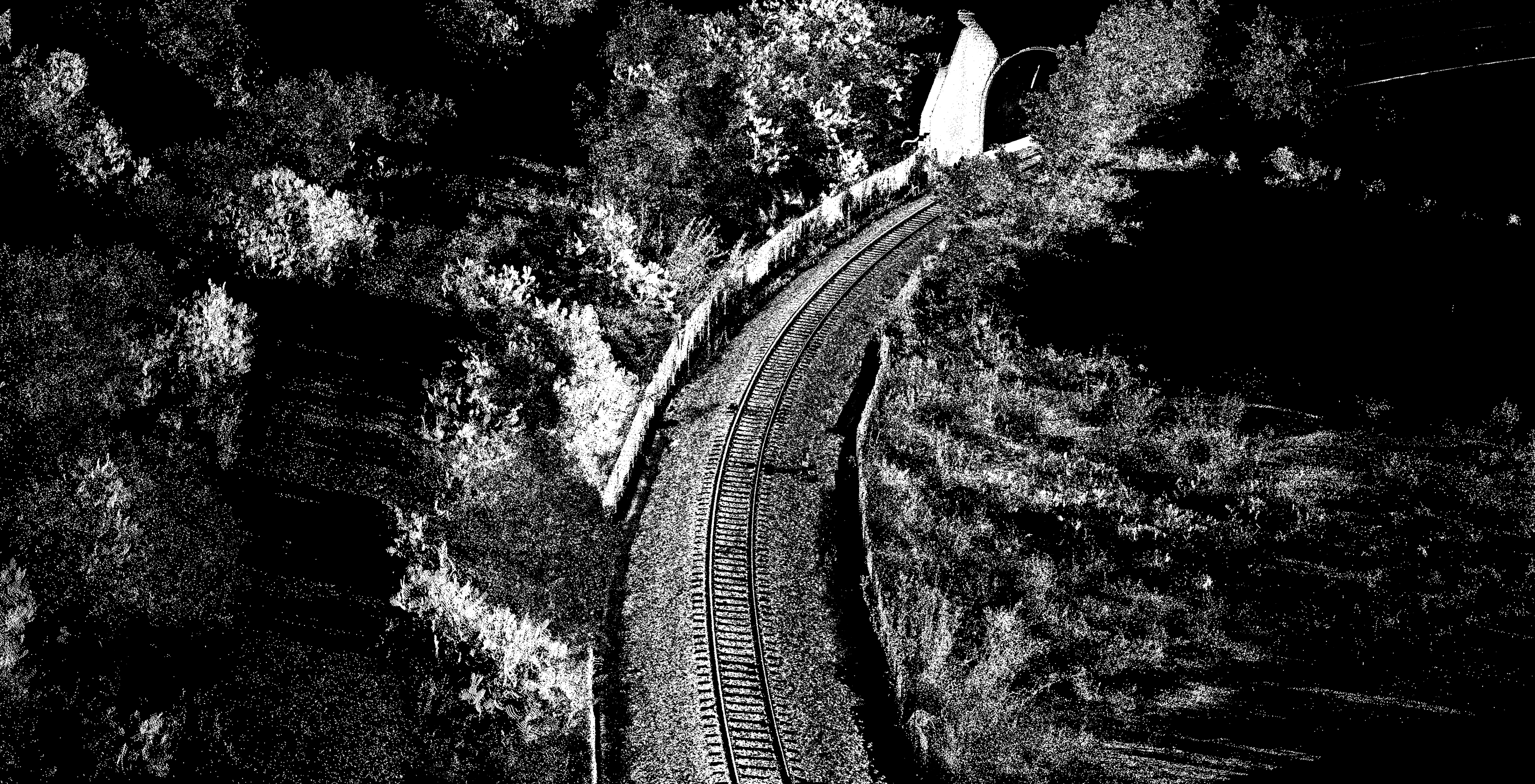}
\end{center}
\vspace*{-.3cm}
\caption{A visual representation of a 3D point cloud showcasing a curved section of a railway emerging from a tunnel, with the surrounding vegetation captured in the scene.}
\label{Fig_rail0}
\end{figure}

In the present example, we delve into a segmentation task performed on a point cloud that portrays a railway environment, captured using a terrestrial laser scanning system. An instance of such a scenario is presented in Figure \ref{Fig_rail0}, which will serve as the subject of our examination. Here, we observe a curved railway emerging from a tunnel, enveloped by dense vegetation. Our aim revolves around identifying the train rails within this scenario and approximating their shape using a suitable spline curve. Conducting such an analysis can yield valuable insights into the transportation system and aid in identifying potential issues that could impact its operational effectiveness (see \cite{ChJuOl19}  and reference therein).

It is worth underscoring that the essence of this example lies in testing the entropy-based approach on a highly noisy dataset, where the set $D_1$ of inliers is significantly dwarfed by the set $D_2$ of outliers. As a result, the technique showcased in this example serves as a proof of concept rather than a definitive solution for the intended problem (for a more effective identification of the rails, refer to works such as \cite{LZTSZC18,Ar15,ADILLLMP22}).

A point cloud is a data set that realizes a digital representation of a physical environment or object in a three-dimensional space. It is arranged in a structured array housing fields that store various attributes for each point within the cloud. These attributes encompass 3D coordinates, distance ranges, color information, intensity measurements, and potentially other geometric or spectral data. We will utilize the intensity parameter, a measure of the reflectivity of the material of the object containing the sample point, to identify reflective elements like train rails.


Within the segmentation procedure, the intensity field frequently comes into play for the purpose of condensing the initial array of data points into a more fitting subset of points pertinent to the analysis. In fact, noteworthy structures, including train rails and overhead wires, exhibit resemblances in their intensity attributes. This correspondence arises from the inherent connection between a surface's reflective characteristics and its constituent material. For instance, train rails are predominantly composed of steel, leading to nearly uniform intensity readings from the laser sensor along the rail's length. By resorting on the intensity parameter as a filtering criterion, we can effectively discern the majority of points situated on the rails.

\begin{figure}[htb]
\begin{center}
\includegraphics[width=0.93\textwidth]{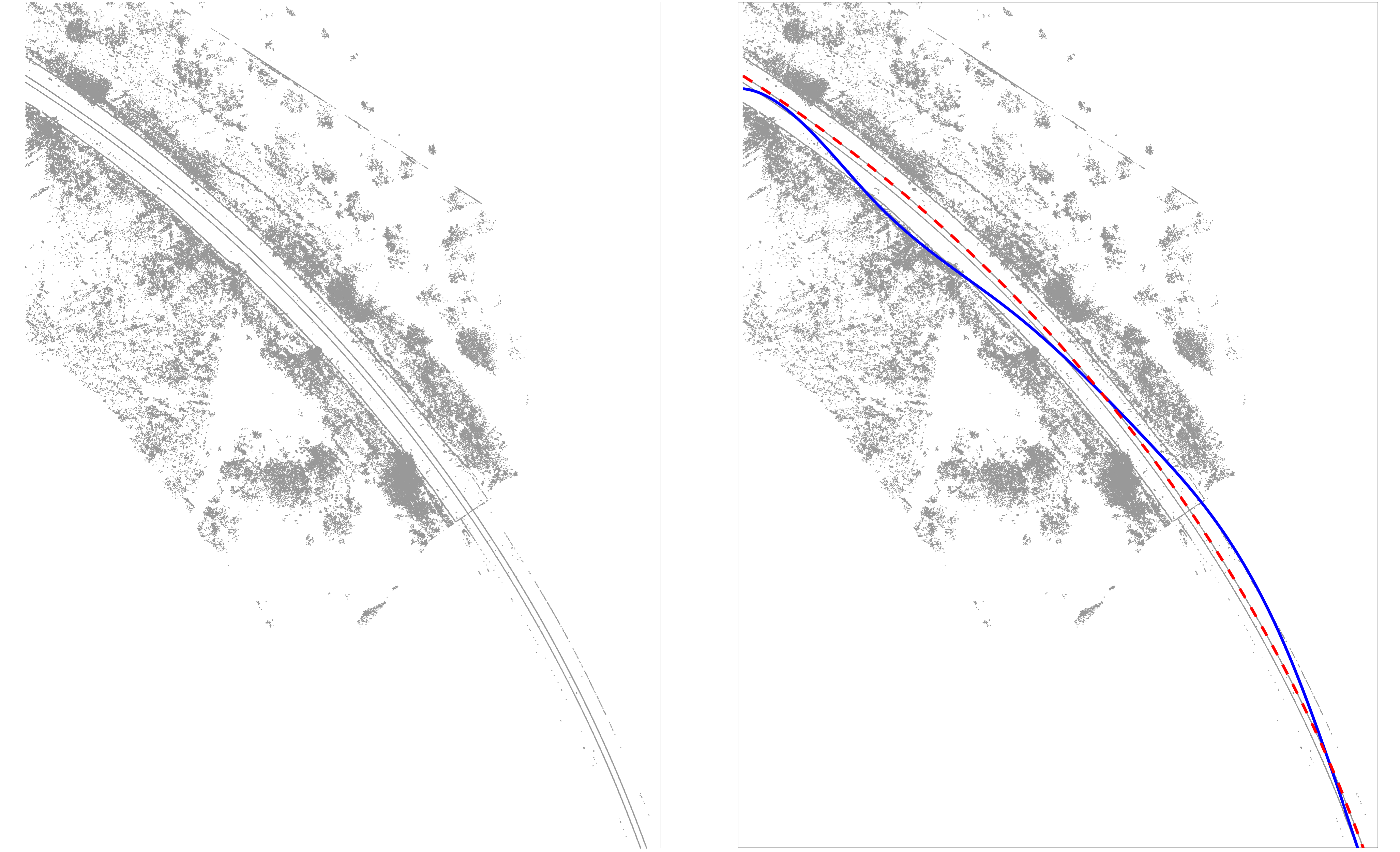}
\vspace*{-.3cm}
\end{center}
\caption{Left picture: 2D projection of the filtered point cloud. The tracks are correctly represented but, unfortunately, vegetation outside the gallery introduces a relevant number of noisy points in the filtered image. Right picture: ordinary least squares spline approximation (solid blue line) and  maximal-entropy least squares spline approximation (dashed red line).}
\label{Fig_rail12}
\end{figure}

Building upon the analysis conducted in \cite{ADILLLMP22} for a point cloud of similar nature, our approach to reduce the size of the original point cloud, while retaining the majority of rail points, involves extracting those with intensity values not exceeding 65. Additionally, due to the level nature of the terrain under consideration, we omit the vertical component of the points and instead focus on a two-dimensional projection of the filtered point cloud. This projection is illustrated in the leftmost image of Figure \ref{Fig_rail12} and forms a data set comprising $304911$ points. The lower-right section of the image corresponds to the segment of the rails situated within the tunnel. This region exhibits a much cleaner appearance compared to the area outside the tunnel. Indeed, in the external environment, a considerable number of points associated with vegetation are regrettably retained even after the filtering procedure. This introduces a notable degree of noise into the data.

The right image in Figure \ref{Fig_rail12} displays the ordinary least squares  spline approximation curve (solid blue line). By referring to equations (\ref{spline})-(\ref{regular}), this curve is obtained through a spline of degree $d=2$ and $n=15$, utilizing a uniform $(d+1)$-regular knots distribution. Evidently, the OLS approximation does not deviate that much from the shape traced by the rail tracks, making it a suitable initial estimate within Algorithm \ref{Alg1} for computing the maximal-entropy weighted least squares spline approximation curve.

This MEWLS curve is depicted in the same graph as a dashed red line. It is clear that the MEWLS spline closely captures the profile of the upper rail, demonstrating a very high accuracy. An inspection of the weights through  formula (\ref{D1D2}) reveals that, in this specific example, the number of outliers exceeds the number of inliers by more than six times.

A comparable process can be subsequently applied to acquire the approximation for the lower rail. This involves eliminating the points related to the upper rail from the dataset and then performing Algorithm \ref{Alg1} again  (we omit to display this latter approximation for visual clarity).

In conclusion, the MEWLS approach effectively enhances the accuracy of the initial OLS approximation and leads to a precise parametric representation of the rails.

\subsection{Detecting and scoring outliers in an environmental data set}
\label{sec:5.3}

The final test case is drawn from a study in \cite{CaCaHo22} and explores an environmental dataset accessible through the R-package \textit{openair} \cite{CaRo12}. This dataset encompasses hourly readings of wind speed, wind direction, and concentrations of pollutants such as NO$_x$, NO$_2$, O$_3$, PM$_{10}$, SO$_2$, CO, and PM$_{25}$ recorded at Marylebone (London) spanning from January 1, 1998, to June 23, 2005. For comparison purposes, we conform to the choice in \cite{CaCaHo22} and  focus on a specific subset of this dataset, only comprising the O$_3$ concentrations during December 2002. This particular segment encompasses a total of 744 observations, while also featuring several instances of missing data points. 

\begin{figure}[htb]
\begin{center}
\includegraphics[width=0.95\textwidth]{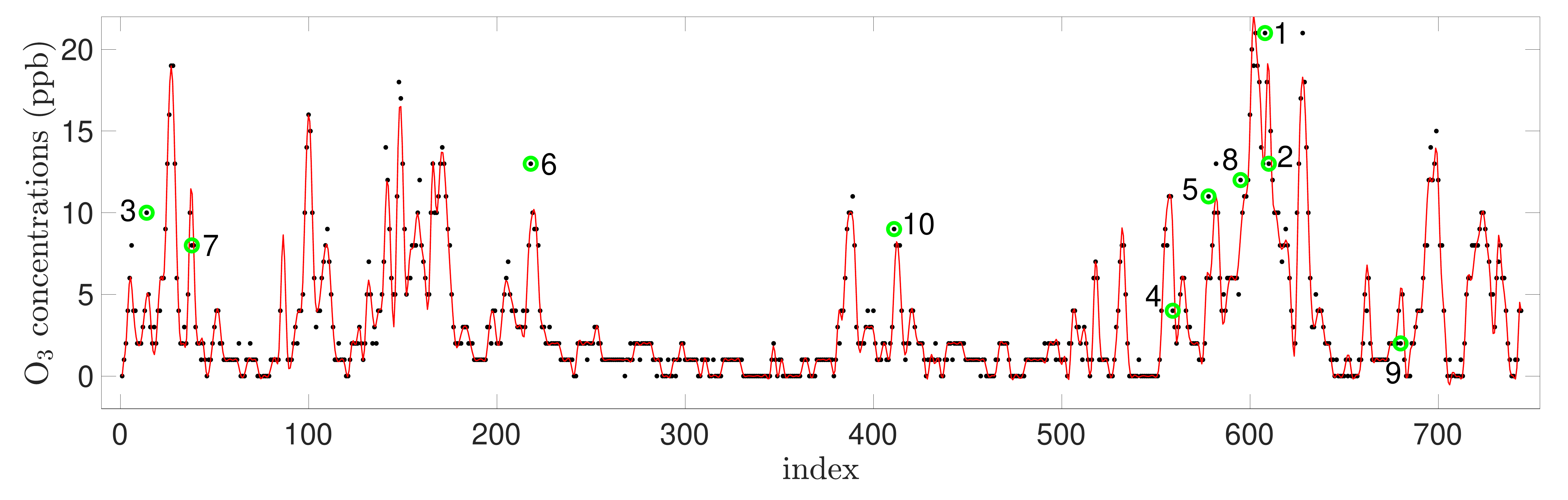}
\end{center}
\vspace*{-.3cm}
\caption{Dots: Hourly O$_3$ concentrations recorded at Marylebone during December 2002 (taken from the R package \textit{openair}). Solid red line: maximal entropy least squares spline approximation. Dots surrounded by green circles identify the first ten outliers detected by the procedure.}
\label{Fig_env}
\end{figure}

The dots depicted in Figure \ref{Fig_env} provide a visual representation of the O$_3$ concentrations, measured in parts per billion (ppb), over the specified time frame. To approximate this univariate time series, we employ a spline function with degree $d=3$, defined on a uniform $(d+1)$-regular knots distribution. In order to capture the erratic nature of the data, we opt for a number of coefficients $n$ equal to half the data points' count. Figure \ref{Fig_env} only displays the MEWLS approximation (red solid line).

In contrast to the approach adopted in prior examples, our strategy for obtaining the approximating spline varies here. Rather than predefining the reduction factor, we pursue a distinct perspective. Specifically, we establish the number of outlier candidates, denoted as $N$, and iteratively reduce the $\MSE$ value until $N$ data points are encompassed within the outlier set $D_2$. This methodology introduces a natural ranking within $D_2$, assigning scores to each prospective outlier. This is readily accomplished using (\ref{D1D2}), where the $i$th point entering $D_2$ receives a score of $i$. In Figure \ref{Fig_env}, outliers are denoted by points enclosed in green circles, each indicating the corresponding score.

The outcomes obtained align with those presented in \cite{CaCaHo22}, particularly those based upon  the extreme value theory. This systematic scoring approach has the potential to streamline the decision-making process, aiding specialists in identifying the data points that merit closer investigation or intervention.


\section{Conclusions}
\label{sec:6}
In real-world scenarios, data quality directly impacts the performance of subsequent analytical processes, so that the importance of effective preprocessing techniques and robust fitting procedures have become increasingly evident.

In this context, we have introduced an entropy-based weighting methodology for determining spline approximations of multivariate time series. In contrast to the ordinary least squares approach, which displays sensitivity to corrupted data, the MEWLS spline approximation effectively mitigates the impact of outliers and noise even when  handling  large and highly noisy datasets. Its ability to accurately extract meaningful information from noisy backgrounds has been illustrated through various synthetic and real-world examples.

One limitation when compared to the OLS approach is that, even for linear models, the resulting algebraic system becomes nonlinear and its solution requires the implementation of an appropriate iterative scheme. In this regard, the OLS solution can serve as an initial estimate.  The numerical illustrations underscore that the MEWLS solution significantly outperforms the classical OLS procedure. Nonetheless, the efficient resolution of this nonlinear system warrants dedicated investigation and will be a focus of future research.

\section*{Acknowledgements}
Felice Iavernaro acknowledges the contribution of the National Recovery and Resilience Plan, Mission 4 Component 2 - Investment 1.4 - NATIONAL CENTER FOR HPC, BIG DATA AND QUANTUM COMPUTING - Spoke 5 - Environmental and Natural Disasters,  under the NRRP MUR program funded by the European Union - NextGenerationEU - (CUP H93C22000450007).

Luigi Brugnano and Felice Iavernaro  thank the GNCS for its valuable support under the INDAM-GNCS project CUP\_E55F22000270001.


%
%


\end{document}